\documentclass[a4paper,11pt,oneside]{article}
\usepackage[top=1.5cm, bottom=1.5cm, left=1.5cm, right=1.5cm]{geometry}

\usepackage{amsmath,amsthm,amsfonts,amssymb,amscd}
\usepackage{subfig}
\usepackage{mathtools}
\usepackage[dvipsnames]{xcolor}
\usepackage{graphicx,bm}
\usepackage{algorithm}
\usepackage{url}
\usepackage{hyperref}
\hypersetup{
    colorlinks,
    linkcolor={red!50!black},
    citecolor={blue!50!black},
    urlcolor={blue!80!black}
}
\usepackage{siunitx}
\usepackage{authblk}
\usepackage[authoryear]{natbib}
\newtheorem{theorem}{Theorem}[section]
\newtheorem{remark}[theorem]{Remark}
\newtheorem{definition}[theorem]{Definition}
\newtheorem{problem}[theorem]{Problem}

\newenvironment{acks}{\paragraph{Acknowledgements}}{}

	\title{ideal.II:\ a Galerkin Space-Time Extension to the Finite Element Library deal.II}
	\author[1,2]{J.P. Thiele}
	
    \affil[1]{Weierstrass Institute Berlin\\
    Mohrenstraße 39, 10117 Berlin, Germany}        
	
    \affil[2]{Institute of Applied Mathematics, Scientific Computing,
    Leibniz University Hannover,
    Welfengarten 1, 30167 Hannover, Germany}       

	\date{}
	
\begin{document}\maketitle
    \begin{abstract}
        The C++ library deal.II provides classes and functions to solve stationary problems with finite elements on one- to threedimensional domains.
        It also supports the typical way to solve time-dependent problems using time-stepping schemes,
        either with an implementation by hand or through the use of external libraries like SUNDIALS.
        A different approach is the usage of finite elements in time as well, which results in space-time finite element schemes.
        The library ideal.II (short for instationary deal.II) aims to extend deal.II to simplify implementations of the second approach.
    \end{abstract}

\section{Introduction}
The idea of using finite elements for the discretization in both
time and space was initially published in~\citep{Od69b}.

Initially, the temporal discretization was done by
constructing a Butcher tableaux, through temporal 
integration by hand, which was then used to apply Rothes method. 
This also allowed for the analysis of these
discretization methods in the context of ODEs
in~\citep{lesaint1974finite,Johnson88,AzizMonk89}.

One of the first direct discretizations by space-time
basis functions was proposed by~\citet{hughes1988space}.
Therein, they applied space-time finite elements
to the solution of linear elastodynamics on a one-dimensional rod
through a native discretization with a two-dimensional unstructured grid.

Recently, with advances in computing power and efficiency,
direct discretizations without the need for deriving
and implementing Butcher tableaus have gained new interest.
Practically, these split into two categories.
\textit{Native discretizations}, as done in~\citep{hughes1988space},
use \(d+1\)-dimensional meshes and finite elements.
These meshes can be fully unstructured, 
i.\,e.\ information about temporal changes
of the spatial domain are embedded.
\textit{Tensor product discretizations} construct meshes
and finite element as products of spatial
and temporal versions.

We will now give a short overview of existing open source software packages
for space-time finite elements. 
As outlined above and detailed in the PhD thesis of the first author~\citep{thiele2024},
there are three major approaches to discretization.
The first, `reformulation as a time-stepping scheme',
is supported by any finite element library as it requires most of the 
temporal discretization to be implemented by hand. 
The second, `native \(d+1\) dimensional dimensional discretizations',
requires support for \(4\)-dimensional meshes and elements for problems 
in three space dimensions.
To the best of the authors knowledge, the only open source finite element library 
providing this at the time of writing is MFEM~\citep{anderson2021mfem}.
The third, `tensor product space-time slabs', is supported by 
any finite element library that offers meshes and elements in \(1\) and \(d\) 
dimensions and we give some more information on using this approach in chapter 3.
For instance, this approach was used in~\citep{hauser2023numerical} using Netgen/NGSolve~\citep{schoberl1997netgen},
in \texttt{Example205\_LowLevelSpaceTimePoisson}
of \texttt{ExtendableFEMBase.jl}\footnote{\url{https://github.com/chmerdon/ExtendableFEMBase.jl}}
and in our own study~\citep{RoThieKoeWi23} using \texttt{deal.II}\citep{dealiiDesignFeatures}.
The code for our own study is openly available\footnote{\url{https://github.com/mathmerizing/dwr-instatfluid}},
was build on previous software by Uwe K\"ocher and can be seen as a predecessor
to \texttt{ideal.II}.

The library \texttt{ideal.II} described in this paper
aims to simplify the assembly and solution of
nonstationary PDEs with the tensor product approach.
It is based on the general purpose library
\texttt{deal.II}~\citep{dealiiDesignFeatures}
and can be seen as an extension to improve usability
in the context of space-time finite elements.

The rest of the paper is structured as follows.
We start by describing the mathematical background
of the tensor product space-time finite element approach.
Then, we will have a short, more general discussion
of how an existing finite element library
can be extended to support tensor product finite elements.
This is followed by more details on the design
and implementation of \texttt{ideal.II}.
Afterwards, we will look at numerical examples that
have been solved using \texttt{ideal.II} to
show some of the functionality it provides.
Finally, we will give an overview of planned features for future versions 
of \texttt{ideal.II} 

\section{Space-time tensor-product finite elements}

\subsection{Basic notation}
Let \(\Omega\subset\mathbb{R}^d\) be the
time-independent spatial domain with dimension \(d\in \{1,2,3\} \)
and Lipschitz-continuous boundary \(\partial\Omega \).
Let \(I=(0,T)\) be the temporal domain with
final point \(T\in\mathbb{R}>0\).
The resulting \textit{space-time cylinder} is denoted by
\(\Sigma = I\times\Omega\subset\mathbb{R}^{d+1}\).
In the following, a bold font, e.\,g. \(\mathbf{f}(x)\in \mathbf{V}(\Omega) \),
will denote a vector-valued function or function space.
\subsection{Bochner spaces}
Using an evolution triple~\citep{zeidler1989linear}
\begin{equation}
    V\subseteq H \subseteq V^*,
\end{equation}
also called Gelfand triple, or rigged Hilbert space \((H,V)\),
we can define the following Bochner spaces
\begin{definition}\label{def:Bochner}
    Let \(V(\Omega)\) and
    \(\mathbf{V}(\Omega)\) be spaces over
    \(\Omega \) which form evolution triples with \(L^2(\Omega)\)
    and \(\mathbf{L}^2(\Omega)\) respectively.
    Then, we can define the function spaces
    \begin{align}
        W(I,V(\Omega))\coloneqq
        \{v\in L^2(I,V(\Omega));
        \partial_t v\in L^2(I,V^*(\Omega))
        \} \\
        \mathbf{W}(I,\mathbf{V}(\Omega))\coloneqq
        \{\mathbf{v}\in L^2(I,\mathbf{V}(\Omega));
        \partial_t \mathbf{v}\in L^2(I,\mathbf{V}^*(\Omega))
        \}.
    \end{align}
    The spaces are both Bochner spaces, i.\,e.
    function spaces with values in Banach spaces~\citep{Ruzicka2020},
    and time-dependent Sobolev spaces~\citep{zeidler1989linear}.
    \begin{theorem}\label{theo:cont_emb}
        For \(W(I,V(\Omega))\) and  \(\mathbf{W}(I,\mathbf{V}(\Omega))\)
        as defined in~\ref{def:Bochner} the following continuous embeddings
        hold
        \begin{align}
            W(I,V(\Omega))\subset C(\bar{I},L^2(\Omega)), \\
            \mathbf{W}(I,V(\Omega))\subset \mathbf{C}(\bar{I},\mathbf{L}^2(\Omega)).
        \end{align}
        \begin{proof}
            See solution to Problem 23.10 in~\citep{zeidler1989linear}.
        \end{proof}
    \end{theorem}
\end{definition}
As an important side note, this result holds not only for
\(H^1(\Omega) = \{ v\in L^2(\Omega); \nabla v\in L^2(\Omega)\} \),
but also for the commonly used spaces
\(\mathbf{H}_\mathrm{curl}(\Omega) = \{ \mathbf{v}\in \mathbf{L}^2(\Omega);
\mathrm{curl}(\mathbf{v})=\nabla\times \mathbf{v}\in \mathbf{L}^2(\Omega)\} \)
and \(\mathbf{H}_\mathrm{div}(\Omega) =
\{ \mathbf{v}\in \mathbf{L}^2(\Omega); \mathrm{div}(\mathbf{v})=\nabla\cdot \mathbf{v}\in L^2(\Omega)\} \).
The latter two function spaces, together with their discrete finite 
element subspaces, play an important part, e.\,g.\ in the solution 
of Maxwell's equations~\citep{Monk2003} or 
weakly divergence free discretizations of Navier-Stokes equations~\citep{JoLiMeNeRe17} respectively, therefore we mention them briefly.

\subsection{Derivation of space-time weak formulations}
With suitable function spaces, we can now derive the space-time weak
formulation, following the same steps as for stationary problems, namely
\begin{enumerate}
    \item Choice of suitable test function spaces \(V(\Omega)\) and \(X=W(I,V(\Omega)) \).
    \item Multiplication of the strong from with test functions \(\varphi\in X \)
          and integration over \(\Sigma \).
    \item If needed, partial integrations in space and/or time.
\end{enumerate}
As an example we will apply these steps to the heat equation
\begin{problem}[Heat equation in strong form]
Find \(u(t,x) \), such that
\begin{equation}
    \begin{aligned}
        \partial_t u(t,x) -\Delta_x u(t,x) & = f(t,x) \text{ in } (0,T)\times\Omega         \\
        u(0,x)                             & = u_0(x) \text{ on } \{0\}\times\Omega         \\
        u(t,x)                             & = g(t,x) \text{ on } (0,T)\times\partial\Omega
    \end{aligned}
\end{equation}
\end{problem}
For step 1 we choose \(V(\Omega) = H^1_0(\Omega)\)
and \(X= X(I,\Omega) = W(I,H^1_0(\Omega))\),
such that step 2 yields

\begin{equation}
    \int\limits_0^T\int\limits_\Omega
    \partial_t u(t,x)\varphi(t,x) \mathrm{d}x\mathrm{d}t
    -\int\limits_0^T\int\limits_\Omega
    \Delta_x u(t,x) \varphi(t,x)\mathrm{d}x\mathrm{d}t
    =\int\limits_0^T\int\limits_\Omega
    f(t,x)\varphi(t,x) \;\forall\varphi\in X(I,\Omega)
\end{equation}
By applying integration-by-parts, i.\,e.\ Greens formula in space on the Laplacian we obtain:
\begin{problem}[Heat equation in continuous weak form]\label{prob:cont_heat}
Find \(u\in X(I,\Omega)+g\), such that
\begin{equation}
    \begin{aligned}
        ((\partial_t u(t,x), \varphi(t,x))) + ((\nabla_x u(t,x),\nabla_x \varphi)(t,x))
               & = ((f(t,x),\varphi(t,x)))\; \forall\varphi\in X(I,\Omega), \\
        u(0,x) & = u_0(x)\text{ on }\{0\}\times\Omega.
    \end{aligned}
\end{equation}
\end{problem}
Here \(((u,v))\) denotes the inner product in \(L^2(I,L^2(\Omega))\)
and in the following \({(u,v)}_{L^2(\Omega)}\) will denote the spatial
inner product, such that
\begin{equation}
    ((u,v)) = \int\limits_0^T {(u,v)}_{L^2(\Omega)}\mathrm{d}t.
\end{equation}
To apply discontinuous Galerkin (dG) methods in time we have to
define a matching function space.
With a given temporal triangulation \(\mathcal{T}_k \) of \(I=(0,T) \) into \(M \) subintervals \(I_m = (t_{m-1},t_m) \),
we allow for discontinuities in the function space and obtain
the \textit{broken Bochner space}
\begin{equation}
    \widetilde{W}(\mathcal{T}_k,V(\Omega))\coloneqq\lbrace u\in L^2(I, L^2(\Omega));\; u|_{I_m}\in W(I_m,V(\Omega))\;\forall I_m\in\mathcal{T}_k \rbrace.
\end{equation}
Theorem~\ref{theo:cont_emb} then yields
\(W(I_m,V(\Omega))\subset C(\bar{I}_m,L^2(\Omega))\),
such that the following limits are well-defined
\begin{align}
    v^+(t_m) & = \lim\limits_{\varepsilon\to 0} v(t_m+\varepsilon)\quad\text{for}\quad t_m\in\lbrace 0,t_1,\ldots , t_{M-1}\rbrace, \\
    v^-(t_m) & = \lim\limits_{\varepsilon\to 0} v(t_m-\varepsilon)\quad\text{for}\quad t_m\in\lbrace t_1,\ldots, t_{M-1},T\rbrace.
\end{align}
With these
we can further define the jump across interval edges as
\begin{equation}
    {[v]}_m\coloneqq v_m^+-v_m^-.
\end{equation}
Using these jump terms and limits from above and below we can weakly enforce the initial condition and derive a discontinuous weak form.

\begin{problem}[Heat equation in discontinuous weak form]\label{prob:discont_heat}
Find \(u\in  \widetilde{X}(\mathcal{T}_k,\Omega)+g \), such that
\begin{equation}
    \begin{aligned}
        \sum\limits_{m=1}^M\int_{I_m}{(\partial_t u, \varphi)}_{L^2(\Omega)} + {(\nabla_x u,\nabla_x \varphi)}_{L^2(\Omega)}\mathrm{d}t
        +\sum\limits_{m=1}^{M-1} {({[u]}_m,\varphi_m^+)}_{L^2(\Omega)} +{(u_0^+,\varphi_0^+)}_{L^2(\Omega)} \\
        = \sum\limits_{m=0}^M\int_{I_m}{(f,\varphi)}_{L^2(\Omega)} +{(u_0,\varphi_0^+)}_{L^2(\Omega)}\; \forall\varphi\in \widetilde{X}(\mathcal{T}_k,\Omega).
    \end{aligned}
\end{equation}
\end{problem}
Note that we also split the integral over \(I \) into integrals over the subintervals due to the discontinuities.
For discontinuous Galerkin methods in space we would conceptually
follow the same steps by introducing a broken Sobolev space
depending on the spatial triangulation, cf.\ \citep{DiPietroErn2011}.

\begin{remark}[Support for cG elements in time]
    Currently, the software is designed to support discontinuous Galerkin discretizations in time.
    Therefore, the discretization and further chapters will focus on dG methods.2
\end{remark}

\subsection{Discretization}\label{subsec:disc}
Due to the structure of \(\Sigma \) as \(I\times\Omega \),
the space \(L^2(I,V(\Omega))\) is isometric to the tensor product space
\(L^2(I)\hat\otimes V(\Omega)\)~\cite[Prop. 1.2.27]{PicardMcGhee2011}.
Here, \(\hat\otimes \) denotes the closure of the Hilbert space tensor product.
Similarly, isometries hold for the (semi-)discrete 
subspaces~\cite[Prop. 1.2.28f]{PicardMcGhee2011}.
With these isometries we can identify space-time basis functions
with the product of temporal and spatial basis functions
\begin{equation}
    \forall\varphi_{tx}\in L^2(I,V(\Omega))\;\exists
    \varphi_t\in L^2(I)\text{ and }\varphi_x\in V(\Omega),
    \text{ such that } \varphi_{tx}(t,x) = \varphi_t(t)\varphi_x(x).
    \label{eq:basis_function_mult}
\end{equation}
Accordingly, application of the product rule yields
\begin{align}
    \partial_t\varphi_{tx}(t,x) & \coloneqq (\partial_t\varphi_t(t))\cdot\varphi_x(x),
    \label{eq:basis_function_mult_dt}                                                  \\
    \nabla_x\varphi_{tx}(t,x)   & \coloneqq \varphi_t(t)\cdot(\nabla_x\varphi_x(x)).
    \label{eq:basis_function_mult_dx}
\end{align}
We can now start discretizing by constructing a finite dimensional
subspace of \(L^2(I)\) using polynomial functions of order \(r\) on each subinterval i.\,e.\  functions in
\begin{equation}
    \mathcal{P}_r(\mathcal{T}_k)\coloneqq
    \lbrace u\in L^2(I);\;u|_{I_m}\in \mathcal{P}_r(I_m)\;\forall I_m\in\mathcal{T}_k\rbrace\subset L^2(I).
\end{equation}
To obtain globally defined functions on \(\Sigma \) we define the temporally discrete function space
\begin{equation}
    \widetilde{X}_k^r \coloneqq L^2(I, L^2(\Omega))\cap (\mathcal{P}_r(\mathcal{T}_k)\otimes H^1_0(\Omega))\subset\widetilde{X}
\end{equation}
To keep the notation consistent with other works, spaces of functions that are allowed to have
discontinuities at the inter-element edges, i.\,e.\  the temporal grid points of \(\mathcal{T}_k\),
are denoted by a tilde.

For spatial functions in \(H^1(\Omega)\) we will follow the standard way
of constructing a conforming finite element space,
which will give us a discretization of the second part of the tensor product.
For the discretization of \(\Omega \) we have two options: fixed or dynamic meshes. \\
On a fixed mesh we use a single spatial triangulation \(\mathcal{T}_h\) 
into (quadrilateral) elements \(K\) and obtain the
space-time triangulation \(\mathcal{T}_{kh}\) by application of a single cartesian product, 
i.\,e.\  \(\mathcal{T}_{kh} = \mathcal{T}_k\times\mathcal{T}_h\).
Then, \(V_h^s(\mathcal{T}_h)\) is the standard space of isoparametric Lagrangian finite elements of order \(s\), i.\,e.\
\begin{equation}
    V_h^s(\mathcal{T}_h)\coloneqq\lbrace v\in V;\; v|_K\in \mathcal{Q}_s(K)\text{ for }K\in\mathcal{T}_h\rbrace \subset H^1_0(\Omega).
    \label{eq:lagrangian_fe}
\end{equation}
Replacing \(H^1_0(\Omega)\) by \(V_h^s(\mathcal{T}_h)\) we obtain the fully discrete function space
\begin{equation}
    \widetilde{X}_{k,h}^{r,s} \coloneqq L^2(I,L^2(\Omega))\cap (\mathcal{P}_r(\mathcal{T}_k)\otimes V_h^s(\mathcal{T}_h))
    \subset\widetilde{X}_k^r.
    \label{eq:fixed_dG_space}
\end{equation}
\begin{remark}
    For use in adaptive mesh refinement in space \(V_h^s\) should allow for hanging nodes when using quadrilateral elements.
\end{remark}
For dynamic meshes we want to be able to have different spatial meshes in time, which we achieve by using one spatial mesh \(\mathcal{T}_h^m\)
per temporal element \(I_m\). For each mesh we obtain a corresponding Lagrangian finite element space \(V_h^{s,m}\).
Then, in~\eqref{eq:fixed_dG_space} the fixed polynomial space  \(V_h^s\)
is replaced by the matching space for each temporal element i.\,e.\
\begin{equation}
    v_{kh}\in \bigcup\limits_{I_m}\left(\mathcal{P}_r(I_m)\otimes V_h^{s,m} \right) \label{eq:polyspace_dynamic}
\end{equation}
\subsection{Resulting space-time shape functions and weak formulations}
For the discontinuous Galerkin method in time, \(cG(s)dG(r)\) in shorthand,
we can now write the fully discrete version of Problem~\ref{prob:discont_heat}
by inserting the appropriate test functions and by replacing the 
continuous boundary term \(g\) by its projection \(\check{g}\) into \(\widetilde{X}_{kh}^{r,s}\)/

\begin{problem}[Fully discrete heat equation in weak form using \(cG(s)dG(r) \) discretization]
Find \(u_{kh}\in \widetilde{X}_{k,h}^{r,s}+\check{g}\), such that
\begin{equation}\label{eq:heat_weak_dg}
    \begin{aligned}
        \sum\limits_{m=1}^M \int\limits_{I_m} {(\partial_t u_{kh}, \varphi_{kh})}_{L^2(\Omega)} + {(\nabla_x u_{kh},\nabla_x\varphi_{kh})}_{L^2(\Omega)}
        \mathrm{d}t
        +\sum\limits_{m=1}^M {({[u_{kh}]}_{m},\varphi_{kh,m}^+)}_{L^2(\Omega)} \\
        + {(u_{kh,0}^+-u_0,\varphi_{kh,0}^+)}_{L^2(\Omega)}
        =  \sum\limits_{m=1}^M \int\limits_{I_m} {(f,\varphi_{kh})}_{L^2(\Omega)}\;\mathrm{d}t \;\forall\varphi_{kh}\in\widetilde{X}_{k,h}^{r,s}.
    \end{aligned}
\end{equation}
\end{problem}

\subsection{Tensor-product slabs}
For the actual solution, the fully discrete weak formulation has to be transformed into a linear equation system.
The approach is mostly the same as for the stationary case with the main difference that the finite element ansatz now is
\begin{equation}
    u_{kh}(t,x) = \sum\limits_{i\in\mathcal{T}_k} \sum\limits_{j\in\mathcal{T}_h} u_{i,j} \varphi_{k,i}(t)\varphi_{h,j}(x)
\end{equation}
with \(i\) and \(j\) as degrees of freedom (DoFs) in time and space respectively.
With some lexicographical ordering \(l(i,j)\) we can instead number all space-time DoFs and obtain
\begin{equation}
    u_{kh}(t,x) = \sum\limits_{l\in\mathcal{T}_{kh}} u_l\varphi_{kh,l}(t,x).
\end{equation}
Using this ansatz we can expand the weak formulation and obtain the
linear equation system with one row and column per space-time DoF.
Solving this single equation system on a fixed mesh
corresponds to a space-time approach that yields the all-at-once solution with a single solve.
However, this needs a lot of memory to store the full matrix
for small discretization parameters \(k\) and \(h\).

As the solution at a temporal DoF \(t\in I_a\) does not depend on the solution
on a future element \(\hat{t}\in I_b\) with \(b>a\), the weak formulation can
be decoupled into single temporal elements \(I_m\).\\
Then, we obtain a set of equation systems \(A_m u_m = F_m\) that have to be solved in sequence. This corresponds to a time-stepping scheme
and could also be rewritten as such.

As a generalization of both time-stepping and all-at-once solutions we introduce tensor product slabs.
These slabs are cartesian products of a set of connected temporal elements with
a single spatial triangulation.
This would allow for space-time multigrid solvers or cheaper construction of preconditioners
compared to the time-stepping approach, while needing less memory compared to the all-at-once approach.
Additionally, in some cases the mesh-to-mesh projection can introduce oscillations.
One example of this is the solution of Stokes or Navier-Stokes problems with Taylor-Hood elements,
where the projection result no longer satisfies the discrete divergence free
condition~\citep{BeWo2012,Linke14}.
Implementation of a slab-based equation system also yields the highest flexibility,
as single element slabs correspond to time-stepping,
while a single slab in time corresponds to a full space-time approach with a fixed spatial mesh.
\section{Extending a 1-to-3d finite element code to tensor product space-time}
\subsection{Prerequisites}
For tensor-product space-time, the code has to provide
the d-dimensional spatial basis functions.
While this paper describes spatial discretization using
finite element methods, the code could in principle also
provide and use a different spatial discretization,
e.g.\, using the finite difference, finite volume, finite cell
or virtual element method.

For the temporal discretization the code has to provide
1-dimensional \(cG\) or \(dG\) basis functions
as well as 1-d meshes.
These two ingredients already allow the solution of non-stationary PDEs
on fixed meshes.

To support dynamic meshes, the library has to provide functions
to interpolate the solution from one mesh to another.

\subsection{Data structures}
In general, we have to differentiate between data structures
for a single slab and for the whole space-time domain.
The latter should be an iterable collection of the former to take
the temporal ordering into account.
\subsubsection{Spacetime structures}
From a convenience standpoint a vector would be a good choice,
but it would store all slab data structures in contiguous memory.
This could potentially result in a lot of reordering, e.\,g.
for storing linear algebra objects like the global solution.
Therefore, a list that only points to the (temporal) successor
is a better choice.
Optimal control methods
and the related DWR method for goal oriented adaptivity
need to solve auxiliary problems backwards in time.
To also allow for a backwards iteration, a doubly-linked list
(\texttt{std::list}) was chosen in \texttt{ideal.II}.

\subsubsection{Slab structures}
Slab data structures fall into two categories.
Firstly, we have objects that wrap existing
\(1\)- and \(d\)-dimensional objects together
and provide a convenient space-time functionality.
Examples of this are triangulation, quadrature rules
and objects for evaluating values space-time finite element
basis functions.
Secondly, there are linear algebra structures which can be
used as-is by using a space-time lexicographical ordering.
Here, we chose what could be termed space-major, i.\,e.
\begin{equation}
    i = i_x + N_x i_t,\label{eq:lexorder}
\end{equation}
with \(N_x \in \mathbb{N}\) and
\(N_t \in \mathbb{N}\) as spatial and temporal
number of DoFs, \(i_x\in [0,\dots,N_x)\) 
and \(i_t\in [0,\dots,N_t)\)  
as spatial and temporal DoF indices and \(i\) as space-time DoF index.\\
This has the advantage of using all the existing linear algebra
functionalities like solvers and preconditioners.
For matrices we also have to think about and prescribe
the correct sparsity pattern \(\mathcal{S}\).
By using~\eqref{eq:lexorder} we can combine a spatial pattern
\(\mathcal{S}_x\) and a temporal pattern \(\mathcal{S}_t\).
Since~\eqref{eq:lexorder} groups all spatial DoF belonging to a
single temporal DoF together, \(\mathcal{S}_t\) reflects
the block structure of the space-time system.

For discontinuous Galerkin finite elements in time,
the off-diagonal entries for the jump terms depend on the chosen support points
or DoF locations.
Typical choices are Gauss-Lobatto, Gauss-Legendre, and Gauss-Radau
quadrature points.

We will look at the resulting temporal
patterns for three \(dG(1)\) elements in time. \\
Gauss-Lobatto (Lo) places both support points on the boundary
of the element, such that the limits coincide with the
entries at the temporal DoFs, resulting in
\begin{equation}
    \mathcal{S}_t^{\text{Lo}}=
    \begin{bmatrix}
        s_{00} & s_{01}                                     \\
        s_{10} & s_{11}                                     \\
               & s_{21} & s_{22} & s_{23}                   \\
               &        & s_{32} & s_{33}                   \\
               &        &        & s_{43} & s_{44} & s_{45} \\
               &        &        &        & s_{54} & s_{55}
    \end{bmatrix}.
\end{equation}
For all other choices at least one limit has to be
constructed by linear combination of all temporal DoFs of the element.
The Gauss-Radau formulae place one DoF on either the left (RL) or right (RR) boundary
and the all other DoFs on the inside of the element, resulting
in the following two patterns
\begin{equation}
    \mathcal{S}_t^{\text{RL}}=
    \begin{bmatrix}
        s_{00} & s_{01}                                     \\
        s_{10} & s_{11}                                     \\
        s_{20} & s_{21} & s_{22} & s_{23}                   \\
               &        & s_{32} & s_{33}                   \\
               &        & s_{42} & s_{43} & s_{44} & s_{45} \\
               &        &        &        & s_{54} & s_{55}
    \end{bmatrix},\;
    \mathcal{S}_t^{\text{RR}}=
    \begin{bmatrix}
        s_{00} & s_{01}                                     \\
        s_{10} & s_{11}                                     \\
               & s_{21} & s_{22} & s_{23}                   \\
               & s_{31} & s_{32} & s_{33}                   \\
               &        &        & s_{43} & s_{44} & s_{45} \\
               &        &        & s_{53} & s_{54} & s_{55}
    \end{bmatrix}.
\end{equation}
Finally, the Gauss-Legendre (Le) formula places both
DoFs on the inside of the element resulting in a full
off-diagonal block.
\begin{equation}
    \mathcal{S}_t^{\text{Le}}=
    \begin{bmatrix}
        s_{00} & s_{01}                                     \\
        s_{10} & s_{11}                                     \\
        s_{20} & s_{21} & s_{22} & s_{23}                   \\
        s_{30} & s_{31} & s_{32} & s_{33}                   \\
               &        & s_{42} & s_{43} & s_{44} & s_{45} \\
               &        & s_{52} & s_{53} & s_{54} & s_{55}
    \end{bmatrix}.
\end{equation}
\subsection{Assembly}
For the assembly of the linear equation system there two options.
The first more manual approach adds a temporal
version of almost everything,
i.\,e.\ quadrature formulae, basis function evaluations
and loops over elements and DoFs.
Only the local contributions to the matrix and right hand side vector are s
space-time.
The weak formulation is then written by using~\eqref{eq:basis_function_mult},~\eqref{eq:basis_function_mult_dt}
and~\eqref{eq:basis_function_mult_dx} explicitly.
The second, more user-friendly, approach directly uses space-time versions
of everything by providing wrapper functions and objects
that handle most details internally,
which is the approach \texttt{ideal.II} was written for.
The following two listings compare the innermost DoF loops
for the assembly of the Laplace term of both approaches
in the case of a single temporal element.
\begin{figure}[H]
    \raggedright{}
    \includegraphics[width=0.75\textwidth]{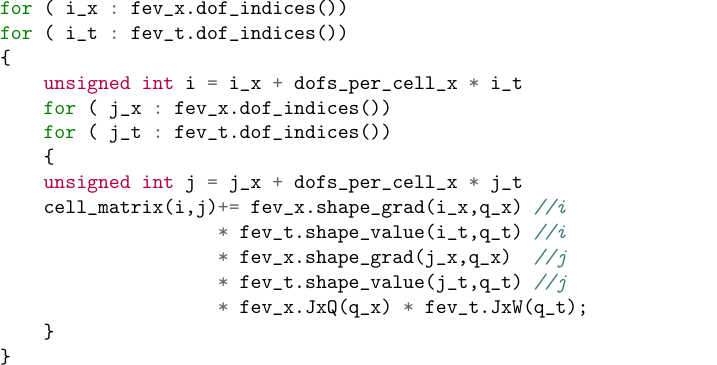}
    \caption{Assembly of \((\nabla u,\nabla\varphi)\)
        in the manual approach using \texttt{deal.II}}\label{list:deal_lap}
\end{figure}

\begin{figure}[H]
    \raggedright{}
    \includegraphics[width=0.75\textwidth]{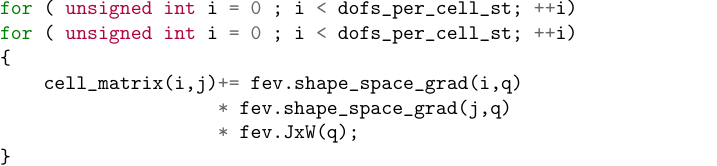}
    \caption{Assembly of \((\nabla u,\nabla\varphi)\) in the
        user-friendly approach using \texttt{ideal.II}}\label{list:ideal_lap}
\end{figure}
\section{Details of ideal.II}
In this section, we start by describing the design decisions
of of \texttt{ideal.II}.
Afterwards, we will look at some implementation details and decisions.
\subsection{Design decisions}
The following decisions and goals went into the design of
the library and it's API.\
\subsubsection{Straightforward space-time weak formulation assembly}
As seen in Listing~\ref{list:deal_lap}, the manual approach
requires a large amount of code duplication.
This makes the code less easy to read and raises the chance
for errors.

To simplify the assembly, the \texttt{FEValues},
\texttt{FEFaceValues} and \texttt{FEJumpValues} classes
offer methods for all space-time functions and derivatives needed.
Internally, the classes call the correct functions of the underlying
\texttt{FEValues} objects of \texttt{deal.II},
which allows usage of all implemented finite element classes in space.
Additionally, the internal code is easily extendable to derivatives
not yet supported as long as they are available in \texttt{deal.II}.

\subsubsection{Final application codes as close as possible
    to stationary \texttt{deal.II} codes}
Unless the system is assembled in an `all-at-once' fashion,
an outer time-marching loop is of course inevitable.
But apart from that, any other functions should be as similar to
a stationary code as possible.
This especially means that the implementation of a specific problem
should not require any by-hand splitting or discretization
of the space-time weak formulation.

The main consequence is that, where possible and sensible,
the number of objects and loops (per slab) stays the same compared
to a stationary code.
As an example, application of Newton's method on a nonlinear
stationary problem typically needs three vectors: the solution,
the Newton update and the residual.
The only additional vector needed for the nonstationary problem is
the one for the initial value of the slab used in the jump term,
which is already needed for linear problems.
For quadrature rules, \texttt{FEValues}, etc.\ there are space-time
versions so the same amount is needed as for a stationary code.
In general this then also applies to loops, as seen in Listing~\ref{list:ideal_lap}.
The only exception are the loops over elements.
As \texttt{deal.II} recalculates local information, like quadrature points, transformation Jacobians,
etc., for the current element, it makes sense to have
an inner loop over temporal elements.
This increases performance by avoiding unneccessary recomputations
of the more expensive spatial elements.
As a side effect this simplifies the application of the
\texttt{WorkStream} concept~\citep{Workstream}
to allow for thread parallel assembly of space-time columns.

\subsubsection{Support for an arbitrary number of intervals per slab}
While we left out the details for multiple temporal
elements in the listings, only a few additional lines of code are needed
in the assembly of nontrivial slabs.
This can be seen in any of the tutorial steps in their
assemble functions.
Namely, the local matrix has to be sized accordingly, the indices have to be offset by the number of
DoFs per space-time element for each temporal element
and the jump term between for inner interfaces between
two temporal elements has to be assembled.
This then makes it easy to change the number of temporal elements
per slab, which may be of interest, e.\,g.\ for specialized preconditioners.
It is possible to limit the maximum number of
temporal elements per slab \(N_{\max} \), which will be considered during temporal refinement.
The edge case of \(N_{\max}=1\) results in a
scheme that corresponds to time stepping,
but with a variable temporal polynomial degree.
The default value of \(N_{\max}=0\) will not split slabs
during refinement, corresponding to an all-at-once approach.

\subsubsection{Extensibility, especially with adjoint problems in mind}
Extensibility is a design goal of \texttt{deal.II} itself,
e.\,g.\ by introducing new types of finite elements.
We aim to leverage and preserve this extensibility because
to us it is one of the most important decisions.
One example of leveraging this is the use of the generic
\texttt{FiniteElement} class for the spatial elements
inside the space-time \texttt{DG\_FiniteElement} class.
This allows users to supply any element (combination)
supported by the underlying version of \texttt{deal.II}
for the spatial discretization.

For optimal control applications or goal oriented adaptivity with the
dual weighted residual method an auxiliary problem has to be solved.
The peculiarity of this problem is that it runs backwards in time.
To allow for the implementation of such adjoint problems,
the underlying data structure for global collection objects
is the \texttt{std::list} as it allows for forwards and backwards iteration
over neighboring entries.
\subsubsection{Preserve most design considerations of \texttt{deal.II}}
The general overview paper~\citep{dealiiDesignFeatures}
outlines the design decisions of \texttt{deal.II} in chapter 2.
Since this library extends \texttt{deal.II} we outline which
of these considerations we will try to preserve in \texttt{ideal.II}.
Additionally, we considered \textit{Design for extensibility}
already since it is one of the most important decisions to us.
\paragraph*{A complete toolbox for finite element codes}
In being tailored for tensor product space-time finite elements,
\texttt{ideal.II} is of course specialized in the use case.
However, for that use case we aim to provide a complete toolbox.
\paragraph*{No hidden magic: deal.II is a library, not a framework}
We aim to provide a library in the same sense.
\paragraph*{Do not reinvent wheels}
Following this rule, many high-level algorithms
and advanced techniques like automatic differentiation
are provided in \texttt{deal.II} by wrapping and using
third party libraries.
Since \texttt{deal.II} is a direct dependency and linked
with all installed third party libraries these algorithms
and techniques can be used in \texttt{ideal.II} as well
and the idea is applied inherently.
Additionally, \texttt{ideal.II} provides helper functions
for calculating DoF indices, sparsity patterns, etc.\
for the \texttt{TrilinosWrappers} interface,
allowing for the implementation of \texttt{MPI}-parallel programs.

\paragraph*{Dimension-independent programming}
This decision is inherited by providing classes
which are templated with the dimension of the spatial
domain \(\Omega \). Internally, they hand this parameter
to the underlying spatial classes of \texttt{deal.II}.
\subsection{Class overview}
Here, we want to give an idea of the general namespace layout of
\texttt{ideal.II} and some details on ideas that went into
the actual implementation of specific classes.

The set of classes for mathematical functions in \texttt{deal.II}
did not need to be wrapped.
This is due to the fact that they have already been designed to work for
time-dependent problems through the provided
\texttt{set\_time()} and \texttt{get\_time()} functions.
These two functions are used internally, e.\,g.\ in\\ \texttt{slab::VectorTools::interpolate\_boundary\_values()}.

One overarching implementation detail is the use of \texttt{std::shared\_ptr}
to store and access references to underlying \texttt{spatial()}
and \texttt{temporal()} classes.
\subsubsection{Namespaces}
The two main namespaces are \texttt{slab} and \texttt{spacetime}.
This differentiates between the classes containing the
actual information for a slab and the collections of
these classes. Examples are \texttt{DoFHandler}
and \texttt{Triangulation}. For vectors there only exist new
\texttt{spacetime} classes, since there is no need for
wrapper classes when using a space-time lexicographical ordering.

In the spirit of this delineation helper functions operating
on \texttt{slab} objects are also located in the \texttt{slab}
namespace. To differentiate their function they are split
further into \texttt{DoFTools} and \texttt{VectorTools} namespaces
following the notation used in \texttt{deal.II}.

Finally, the \texttt{spacetime} namespace also contains
classes which don't need local information during their creation
and can therefore be seen as invariant to a specific slab.
They can of course need and use local information during
their operation.
Examples are quadrature formulae, \texttt{DG\_FiniteElement}
and the family of \texttt{FEValues} classes.

\subsubsection{Triangulation}
In~\ref{subsec:disc} we have seen that there are
two options for constructing a space-time triangulation,
fixed and dynamic.
Currently, only the fixed triangulations are implemented,
but a virtual base class provides a common interface
for future extensions.
For the fixed triangulation the slab objects are constructed
by sharing a pointer to a single common spatial \texttt{dealii::Triangulation}
object. This significantly reduces the memory requirement on fixed
meshes.
Additionally, there are versions with \texttt{MPI}-parallel spatial
triangulations
based on \texttt{dealii::parallel::distributed::Triangulation}.
\subsubsection{DG\_FiniteElement}
This class combines the references to a spatial
and temporal \texttt{FiniteElement} into the
space-time basis functions on the reference element \({(0,1)}^{d+1}\).
The user has full control over the spatial element
since the constructor expects a \texttt{std::shared\_ptr} to
an object of the base class \texttt{FiniteElement}, i.\,e.\ any element
including \texttt{FESystem} for coupled problems.

The temporal basis functions are constructed inside
the \texttt{DG\_FiniteElement} class, depending on the given
temporal degree \(r\) and the \texttt{support\_type},
which can either be \texttt{Lobatto} (default),
\texttt{Legendre}, \texttt{RadauLeft} or \texttt{RadauRight}.
\subsubsection{FEValues family}
In terms of simplifying the direct assembly of
space-time weak formulations, this family of classes
is the centerpiece that does the heavy lifting
and manages the local information needed on each element.
The constructor looks like the \texttt{deal.II} counterpart,
but expects \texttt{spacetime::} versions apart from the
\texttt{dealii::UpdateFlags} passed to both internal objects.

For performance reasons the calculation of
space-time element specific values is split such that\\
\texttt{reinit\_space()} and \texttt{reinit\_time()}
get a reference to the matching sub-element and pass it
to the respective underlying \texttt{dealii::FEValues} object.

\paragraph*{Local indices and global mapping}
The local lexicographical ordering uses the number
of space-time DoFs \(n_{xt}\) per element instead of the total number.
Otherwise, it stays the same, such that the DoF index
shifted by multiples of \(n_{x}\) belongs to the same spatial DoF
at a later temporal DoF.
The easiest way to handle local linear algebra
is the construction of a full matrix with a size
equal to the number of space-time Dofs per element
times the number of temporal elements.
However, for a large slab this matrix will have a lot of
zero entry blocks and could be optimized.

For the global assembly, i.\,e.\ the addition of local
values into the global matrices and vectors,
the function
\texttt{get\_local\_dof\_indices(std::vector<dealii::types::global\_dof\_index>\&)}
returns the mapping from local indices with the above ordering
to the global indices.

\paragraph*{Function evaluation and quadrature }
As explained above, time-dependent mathematical functions are already
provided by \texttt{deal.II} through \texttt{set\_time(double t)}.
The needed temporal quadrature point location \(t\)
is obtained using
the method \texttt{time\_quadrature\_point()}, which
expects a space-time quadrature index.
Similarly, \texttt{space\_quadrature\_point()} returns
the spatial \texttt{dealii::Point} at the given index.

To assemble the local matrices and vectors
we also need \texttt{JxW()}, which combines the quadrature weights and Jacobian determinant of the
transformation from the reference element to the actual element
at a given space-time quadrature point.

\paragraph*{Evaluation of space-time basis functions}
The final component needed for the assembly of element local
objects is the evaluation of space-time basis functions
and function derivatives at each quadrature point.

For scalar problems the following functions are currently
implemented in \texttt{FEValues}
\begin{itemize}
    \item \(\varphi(t,x)\) as
          \texttt{double shape\_value(i,q)}
    \item \(\partial_t\varphi(t,x)\) as
          \texttt{double shape\_dt(i,q)}
    \item \(\nabla\varphi(t,x)\) as
          \texttt{dealii::Tensor<1,dim> shape\_space\_grad(i,q)}
\end{itemize}
with space-time DoF and quadrature indices \(i\)
and \(q\).

Additionally, we need functions to calculate the limits \(\varphi^+\)
and \(\varphi^-\) for each spatial DoF of the given element.
They are provided by the \texttt{FEJumpValues} class through
\texttt{shape\_value\_plus()} and\\ \texttt{shape\_value\_minus()}.

\paragraph*{Evaluation of finite element functions}
When the evaluation of a given finite element function \(u_{kh}\),
e.\,g.\ for nonlinear problems is needed the
following functions evaluate the finite element
ansatz \( u_{kh}(q) = \sum\limits_{i\in K} u_{kh,i}*
\varphi_i(q)\) for each quadrature point \(q\) on the given space-time element \(K\).
\begin{itemize}
    \item
          \texttt{get\_function\_values(ukh, std::vector<double>\&)}
    \item
          \texttt{get\_function\_dt(ukh, std::vector<double>\&)}
    \item
          \texttt{get\_function\_space\_gradients(ukh, std::vector<dealii::Vector<double>>>\&)}
\end{itemize}
with \texttt{dealii::Vector<double> ukh} as finite element function.

\paragraph*{Coupled and vector valued problems}
As mentioned in the design goals, the spatial component of a
\texttt{spacetime::DG\_FiniteElement} can be an arbitrary
\texttt{dealii::FiniteElement} including
\texttt{FESystem}.
In the following, we will look at the example of a
Taylor-Hood \(\mathbf{Q}_2/Q_1\) element
that is typically used for solving (Navier-)Stokes equations 
with a vector valued velocity \(\mathbf{v}\in\mathbf{Q}_2\)
and scalar valued pressure \(p\in Q_1\). For the full
problem statement and code we refer to the tutorial step
22 of \texttt{deal.II} for the stationary case and
step 2 and 4 of \texttt{ideal.II} for the nonstationary case.

To specify the component of the \texttt{FESystem}
which we want to evaluate \texttt{deal.II} provides the classes
\texttt{FEValuesExtractors}, e.\,g.\
\texttt{Scalar} and \texttt{Vector}.
During construction they are given the first index of the component
in the \texttt{FESystem}, e.\,g\
\texttt{FEValuesExtractors::Vector velocities(0)}.

In \texttt{deal.II} the \texttt{[]} operator of
\texttt{FEValues} expects an extractor and
returns a reference to the matching \texttt{FEValuesViews} object.
These offer specialized version of common evaluations,
e.\,g.\ \texttt{value()} or \texttt{gradient()},
with return types of matching dimension.
They also offer further applicable differential operators
like \texttt{divergence()} or \texttt{curl()}.
Internally, this works by constructing all possible
combinations of \texttt{FEValuesViews}
into an additional cache during \texttt{reinit()}.

To provide this syntax in \texttt{ideal.II}, another cache
would be needed, but it would mostly store redundant information
and produce unneccessary memory requirements.
It would also produce an internal implementation that is not easy
to read and extend when needed.

For these reasons we decided to change the syntax
into a single function call where the extractor is the
first argument, e.\,g.\
\texttt{fe\_values.vector\_divergence(velocities,i,q)}.
Instead of the general \texttt{[]} operator handling
the specification the functions are now specialized
by a matching prefix.

\subsubsection{TimeIteratorCollection}
To increase readability of the time marching loop
and to avoid missing part of an iterator over one of the spacetime
\texttt{std::list} collections, \texttt{ideal.II} provides
the \texttt{TimeIteratorCollection}.
Any of these iterators and their underlying collection can be registered with the
TIC through \texttt{add\_iterator()}
and all registered iterators
will be incremented or decremented together with the matching
functions. To check if any of the iterators is at one of
the ends \texttt{at\_end()} and \texttt{before\_begin()}
are provided and serve as stopping criteria for the time marching \texttt{for} loops.

\section{Numerical examples using ideal.II}
With these examples we want to show the current capabilities of the software by discussing the results of two tutorial programs.

Furthermore, we want to show how different temporal and spatial discretization order behave under uniform mesh refinement.

For the discretization we will vary the orders \(r\) and \(s\).
And for the first example we will also vary the type of support points
for the temporal refinement study.
Additionally, we did see benefits from using Gauss-Legendre support
points in a previous study~\citep{RoThieKoeWi23}, which was our main motivation
for supporting these in \texttt{ideal.II}

\subsection{Heat equation}
In this first example we will solve the heat equation~\eqref{eq:heat_weak_dg}.
The problem setting described below is solved in tutorial step 3 of \texttt{ideal.II},
with additional scripts for automation of parameter sweeps and for producing the plots.
These scripts together with the output results can be found in~\citep{idealii-step-3}.
\subsubsection{Problem statement}
This example uses a manufactured solution designed in~\citep{hartmann1998posteriori} for testing adaptive mesh refinement on
dynamical meshes. In previous work we studied this configuration using goal-oriented adaptive 
refinement~\citep{ThiWi24}.
While ideal.II was designed with this application in mind, the initial version so far only supports
fixed meshes.
The solution reads as
\begin{equation}
    \label{eq:hartmann_solution}
    u(t,x,y) = \frac{1}{1+50({(x-\frac12-\frac14\cos(2\pi t))}^2+{(y-\frac12-\frac14\sin(2\pi t))}^2)}
\end{equation}
and is well suited for testing various temporal orders as it is nonlinear in time.
\subsubsection{Configuration}
The equation is solved on the unit cube space-time cylinder, i.\,e.\  \(\Omega = {(0,1)}^2\) and \(T=1\) and initial and boundary condition
follow directly from the given solution~\eqref{eq:hartmann_solution}.
\subsubsection{Quantities/Metrics}
As the exact solution is known, we are calculating the space-time \(L^2\) error i.\,e.\
\begin{equation}
    ||u||_{L^2(I\times\Omega)} = {\left(\int_I ||u(t)||_{L^2(\Omega)}^2\;\mathrm{d}t \right)}^{1/2}
\end{equation}
The error can be split into a temporal and a spatial part which are expected to be
independent of refinement of the other mesh component, such that temporal refinement does not
influence the spatial error and vice versa.
Therefore, we try to choose a suitably small initial \(k\) or \(h\) for the \(h\)- or \(k\)-refinement
studies.
\subsubsection{Results}
Since the resulting number of DoFs on the same mesh differ between polynomial orders,
the following results will not always use the same initial spatial and/or temporal mesh.
Instead, the meshes were chosen in a way that the number of unknowns is similar enough to
be well represented in a single plot. For the comparison of the convergence order
we have done all simulations with Gauss-Lobatto support points. 
To compare the impact of this choice, 
we have also run the temporal refinement with the 
other three choices.

\begin{figure}[H]
    \begin{center}
        \includegraphics[page =1,width=0.8\textwidth]{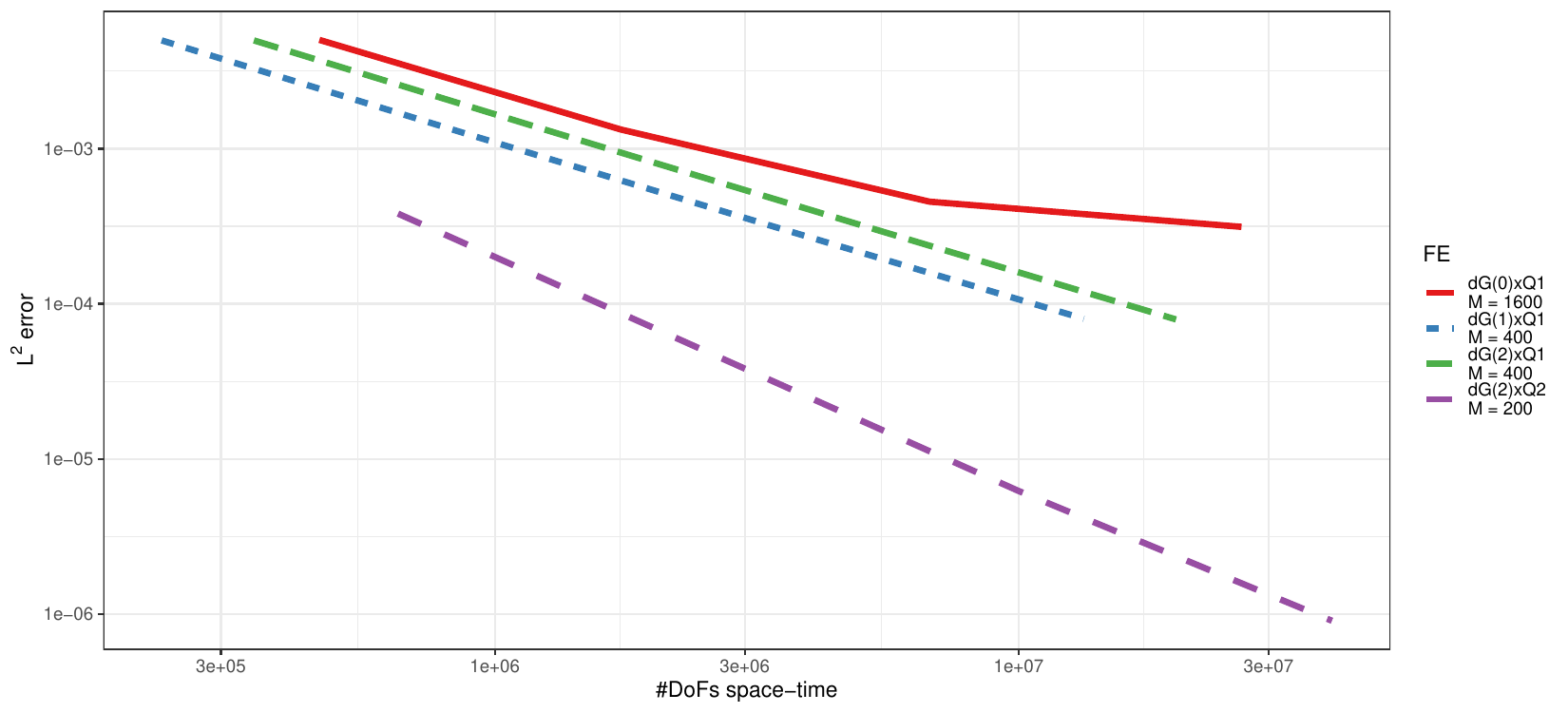}
    \end{center}
    \caption{Results for uniform \(h\)-Refinement i.\,e.\  only in space}\label{fig:plot_hartmann_h_ref}
\end{figure}

Figure~\ref{fig:plot_hartmann_h_ref} shows the results under uniform spatial refinement.
As expected, the convergence orders for bilinear elements in space are identical,
resulting in parallel lines
and the biquadratic elements show a higher order.
However, due to the bad temporal approximation
with piecewise linear elements, the curve for \(r=0\) gets dominated by the
temporal error part.
To remedy this, a much smaller initial \(k\) would be needed,
but then the respective curve would be outside of the 
current plot limits and harder to compare.

\begin{figure}[H]
    \begin{center}
        \includegraphics[page =2,width=0.8\textwidth]{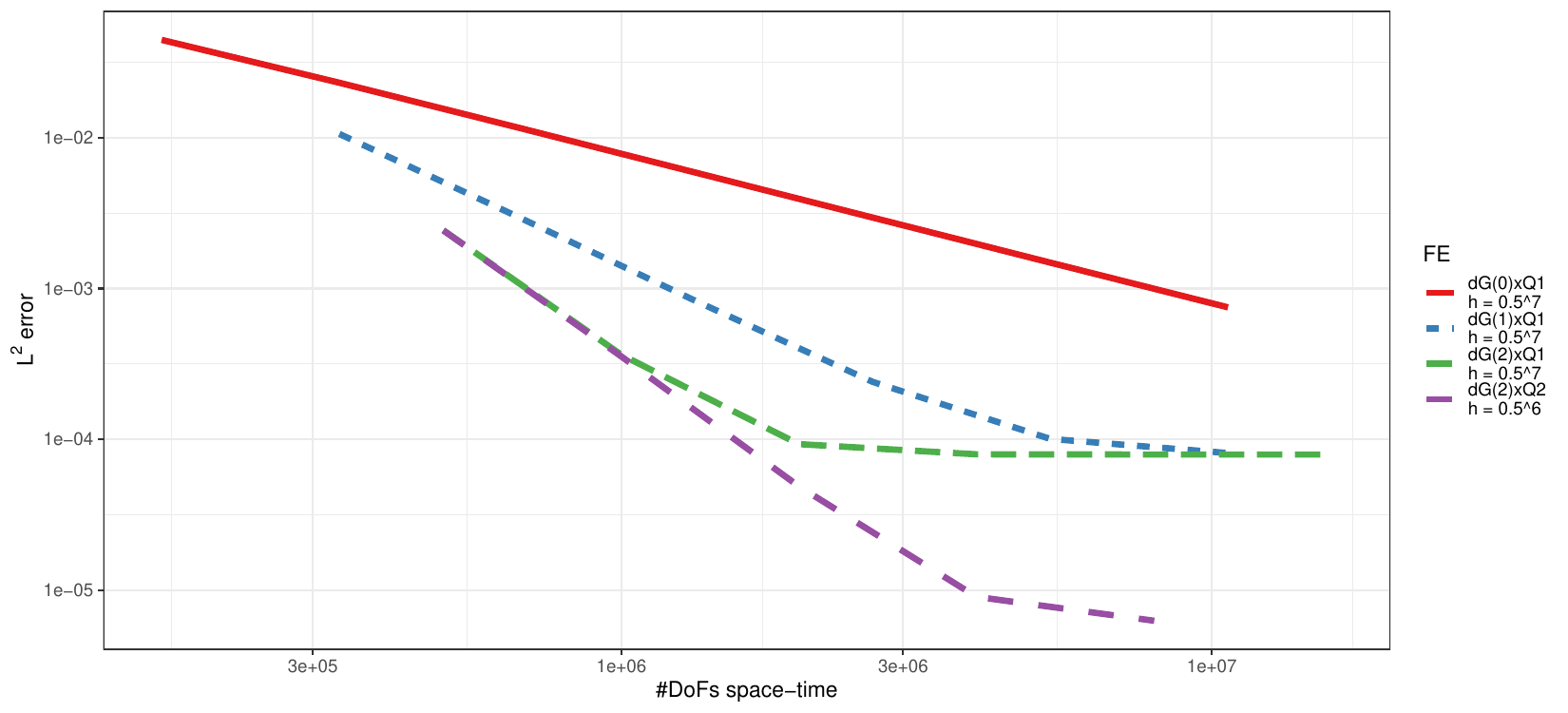}
    \end{center}
    \caption{Results for uniform \(k\)-Refinement i.\,e.\  only in time}\label{fig:plot_hartmann_k_ref}
\end{figure}

The results for uniform temporal refinement are shown in Figure~\ref{fig:plot_hartmann_k_ref}.
There, we see that a higher element order in time also leads to a better convergence rate.
Additionally, we see that both curves with bilinear elements in space
and higher order elements in time converge to the same overall error,
which corresponds to the independence of the spatial error on temporal refinement.

\begin{figure}[H]
    \begin{center}
        \includegraphics[page =3,width=0.8\textwidth]{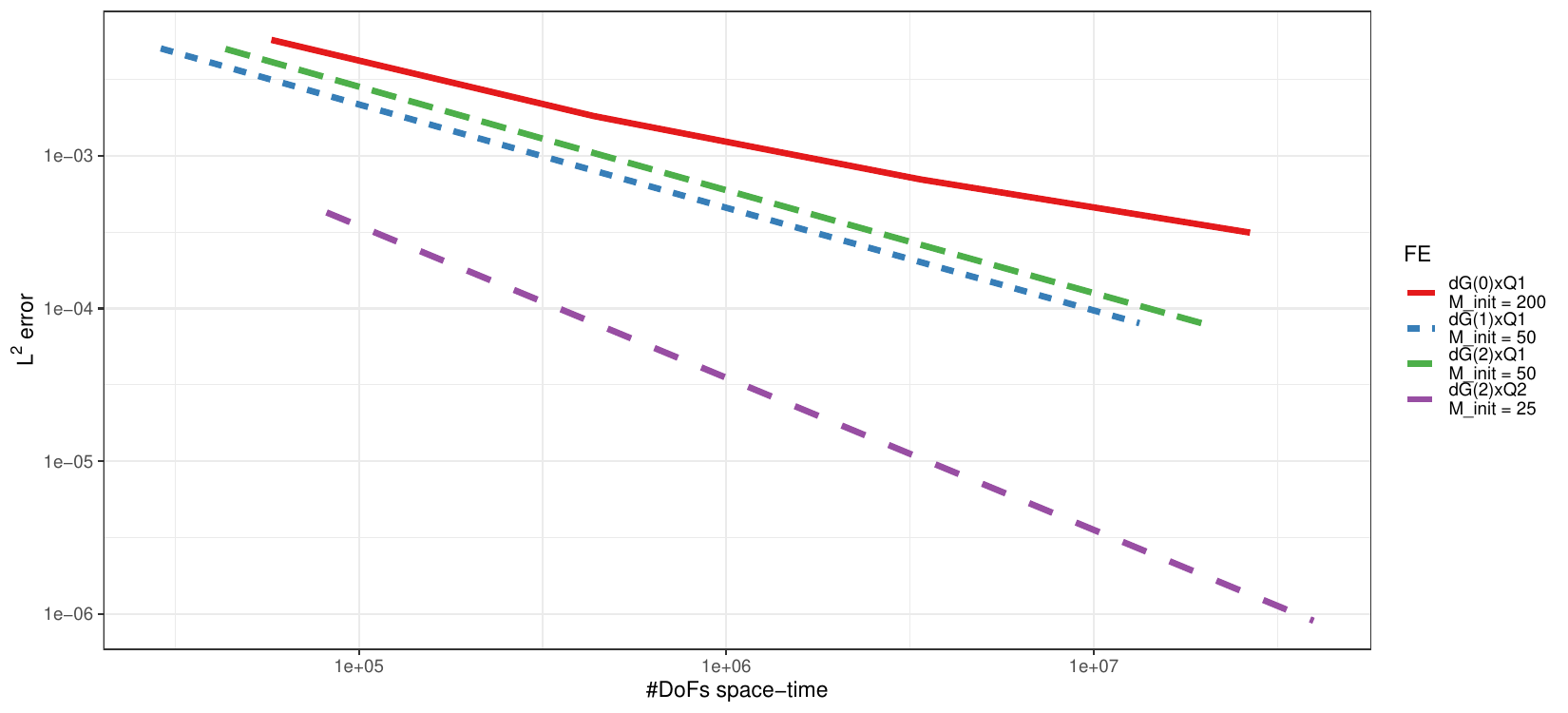}
    \end{center}
    \caption{Results for uniform \(kh\)-Refinement i.\,e.\  in space and time simultaneously}\label{fig:plot_hartmann_kh_ref}
\end{figure}

Combining the previous findings, we see the results for simultaneous uniform refinement
in Figure~\ref{fig:plot_hartmann_kh_ref}.
Here, the overall error is limited by lower order elements,
such that biquadratic elements in space only yield a better
convergence rate when matching quadratic elements in time are chosen.

Finally, we want to compare the impact that the choice 
of support point has on the \(L^2(I,L^2(\Omega))-\)error.
For this we have run all uniform temporal refinements 
for \(r>0\) with Gauss-Legendre and both Gauss-Radau support points 
and calculated the quotient each run and the respective errors obtained from 
the previous simulations with Gauss-Lobatto.
Tables~\ref{table:heat_k_stypes_cg1dg1} to~\ref{table:heat_k_stypes_cg2dg2}
show the proportion between the errors.
A smaller proportion means that the obtained error was smaller then the 
one obtained with Gauss-Lobatto support points.
Overall we can see that the difference between the different 
support type choices gets smaller when using higher order 
elements or a finer mesh and that the overall convergence 
order is not affected.
Comparing the different types on a fixed refinement level, 
we can see that the left Gauss-Radau points consistently yield the 
lowest errors, followed by Gauss-Legendre and the right 
Gauss-Radau rule. We can also see that each of these choices
performs better than Gauss-Lobatto.
Overall, the differences are not huge but we have to keep 
in mind that this is still a relatively simple PDE and problem 
to solve. The differences and therefore advantages in using 
other support points than Gauss-Lobatto may be larger 
for more complicated problems.

\begin{table}[H]
    \centering
    \caption{Resulting error proportions between the chosen support types 
    and Gauss-Lobatto after k-refinment with \(cG(1)dG(1)\) elements.}\label{table:heat_k_stypes_cg1dg1}
    \begin{tabular}{rccc}
     \#DoFs & 
     Left Gauss-Radau & 
     Gauss-Legendre & 
     Right Gauss-Radau \\ \rule{0pt}{14pt}
       \(    332\,820\) & \( 96.62\% \) & \( 97.06\% \) & \( 98.52\% \) \\
       \(    665\,640\) & \( 96.86\% \) & \( 97.62\% \) & \( 99.03\% \) \\
       \( 1\,331\,280\) & \( 97.37\% \) & \( 97.98\% \) & \( 99.02\% \) \\
       \( 2\,662\,560\) & \( 97.87\% \) & \( 98.20\% \) & \( 98.77\% \) \\
    \end{tabular}
\end{table}

\begin{table}[H]
    \centering
    \caption{Resulting error proportions between the chosen support types 
    and Gauss-Lobatto after k-refinment with \(cG(1)dG(2)\) elements.}\label{table:heat_k_stypes_cg1dg2}
    \begin{tabular}{rccc}
     \#DoFs & 
     Left Gauss-Radau & 
     Gauss-Legendre & 
     Right Gauss-Radau \\ \rule{0pt}{14pt}
       \(    499\,230\) & \( 99.22\% \) & \( 99.30\% \) & \( 99.85\% \) \\
       \(    998\,460\) & \( 99.25\% \) & \( 99.31\% \) & \( 99.71\% \) \\
       \( 1\,996\,920\) & \( 99.79\% \) & \( 99.81\% \) & \( 99.91\% \) \\
       \( 3\,993\,840\) & \( 99.75\% \) & \( 99.79\% \) & \( 99.99\% \) \\
    \end{tabular}
\end{table}

\begin{table}[H]
    \centering
    \caption{Resulting error proportions between the chosen support types 
    and Gauss-Lobatto after k-refinment with \(cG(2)dG(2)\) elements.}\label{table:heat_k_stypes_cg2dg2}
    \begin{tabular}{rccc}
     \#DoFs & 
     Left Gauss-Radau & 
     Gauss-Legendre & 
     Right Gauss-Radau \\ \rule{0pt}{14pt}
       \(    499\,230\) & \( 99.22\% \) & \( 99.30\% \) & \( 99.86\% \) \\
       \(    998\,460\) & \( 99.23\% \) & \( 99.29\% \) & \( 99.72\% \) \\
       \( 1\,996\,920\) & \( 99.32\% \) & \( 99.42\% \) & \( 99.72\% \) \\
       \( 3\,993\,840\) & \( 99.75\% \) & \( 99.80\% \) & \( 99.92\% \) \\
    \end{tabular}
\end{table}

\subsection{Navier-Stokes equations}
In this second example we solve the nonlinear Navier-Stokes equations
given below. The code to solve this example is given in
tutorial step 4. It can be found together with the accompanying
scripts in~\citep{idealii-step-4}.

\begin{problem}[Navier-Stokes equations in strong form]
Find vector valued velocity \(\mathbf{v}\) and scalar pressure \(p\),
such that:
\begin{equation}
    \begin{aligned}
        \partial_t\mathbf{v} - \nabla\cdot(\nu\nabla\mathbf{v})
        +\nabla p+(\mathbf{v}\cdot)\mathbf{v} & = 0 \\
        \nabla\cdot v                         & = 0
    \end{aligned}
\end{equation}
\end{problem}
We will look at the specific problem of a channel flow,
for which typical boundary conditions are 
\begin{itemize}
    \item \(\mathbf{v} = 0\) (no-slip) on the channel wall 
    \(\Gamma_\text{wall}\)
    \item \(\mathbf{v} = \mathbf{v}_\text{in}\) at the inflow boundary 
    \(\Gamma_\text{in}\)
    \item \(\partial_n\mathbf{v} = 0\) on the outflow boundary 
    \(\Gamma_\text{out}\) 
\end{itemize}

For the derivation of the weak formulation we use the spatial
space \(\mathbf{V}(\Omega) = \{\mathbf{v}\in \mathbf{H}^1_0(\Omega); p\in L^2(\Omega)\} \).
Discretization in space is done with inf-sup stable Taylor-Hood
element which is a mixed space of two Langrangian
finite elements~\eqref{eq:lagrangian_fe}, i.\,e.
\(\mathbf{v}\in\mathbf{V}_h^{2},p\in V_h^{1}\),
as described in~\citep{GiRa86}.
\begin{problem}[Navier-Stokes equations in weak form]
Find \(\mathbf{u}_{kh}={[\mathbf{v}_{kh},p_{kh}]}^T\in (\widetilde{\mathbf{X}}_{kh}^{r,s_{\mathbf{v}}}+\vec{g}_D^{\vec{v}})\times(\widetilde{X}_{kh}^{r,s_p})\), such that
\begin{equation}
    \label{eq:weak_nse_dG}
    \begin{aligned}
        A_\text{NSE}^{dG}(\mathbf{u}_{kh},\mathbf{w}_{kh}) \coloneqq & \sum\limits_{m=1}^M
        \int_{I_m}
        (\partial_t \mathbf{v}_{kh},\varphi_{kh}^{\mathbf{v}})_{L^2(\Omega)} 
        +\nu(\nabla \mathbf{v}_{kh},\nabla\varphi_{kh}^{\mathbf{v}})_{L^2(\Omega)} 
        +(\mathbf{v}_{kh}\cdot \nabla \mathbf{v}_{kh},\varphi_{kh}^{\mathbf{v}})_{L^2(\Omega)} 
        \\
        & -(p_{kh},\nabla\varphi_{kh}^{\mathbf{v}})_{L^2(\Omega)} 
        +(\nabla\cdot\mathbf{v}_{kh},\varphi_{kh}^p)_{L^2(\Omega)} 
        \;\mathrm{d}t                                                                                                                                             
        \\
        & +\sum\limits_{m=1}^{M-1} ([\mathbf{v}_{kh}]_m,\varphi_{kh,m}^{\mathbf{v},+})_{L^2(\Omega)} 
        +(\mathbf{v}_{kh,0}^+,\varphi_{kh,0}^{\mathbf{v},+})_{L^2(\Omega)} \\ 
        =       
        & (\mathbf{v}^0,\varphi_{kh,0}^{\mathbf{v},+})_{L^2(\Omega)}                                 
        \eqqcolon F_\text{NSE}^{dG}(\vec w_{kh})\;\forall \mathbf{w}_{kh}\in \widetilde{\mathbf{X}}_{kh}^{r,s_{\mathbf{v}}}\times\widetilde{X}_{kh}^{r,s_p},
    \end{aligned}
\end{equation}
with \(\mathbf{w}_{kh}=[\varphi_{kh}^{\mathbf{v}},\varphi_{kh}^p]^T\). 
\end{problem}
\subsubsection{Problem statement}
\begin{figure}[H]
    \begin{center}
        \includegraphics[width=0.95\textwidth]{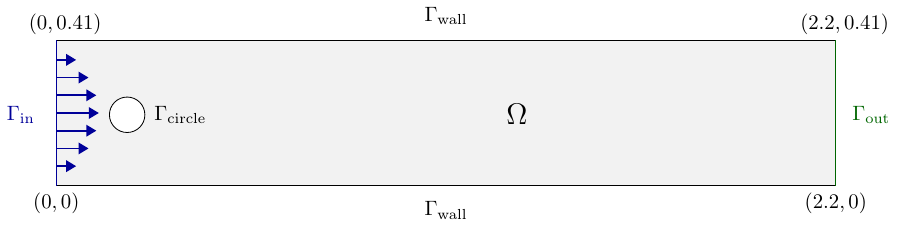}
    \end{center}
    \caption{Spatial domain of the 2D-3 flow around a cylinder benchmark}\label{fig:channel}
\end{figure}

The flow around a cylinder 2D-3 benchmark is 
one of a set of problems introduced and studied 
in~\citep{TurSchabenchmark1996}.
All 2D benchmarks therein describe a Poisseuile flow in a channel of 
height \(2.2\) and length \(0.41\) with a cylindrical 
obstacle of diameter \(0.1\) placed at \((0.2,0.2)\) as shown in Figure~\ref{fig:channel}.
In the 2D-3 case the inflow profile (at \(\Gamma_\text{in}\)) is scaled in time
by the positive half of a sine wave over the interval \(I=(0,8)\).
The time dependent inflow profile is given as
\begin{equation}
\mathbf{v} = \begin{pmatrix}
\sin(\pi t/8)(6y(0.41-y))/(0.41^2)\\
0
\end{pmatrix}
\end{equation}
The walls (\(\Gamma_\text{wall}\) and \(\Gamma_\text{circle}\)) have a no-slip boundary condition \(\mathbf{v}=0\)
and the outflow (\(\Gamma_\text{out}\)) has a natural boundary condition \(\partial_\mathbf{n}\mathbf{v}=0\)

\subsubsection{Quantities/Metrics}
As stated in the original paper we are interested in the time curve of the drag and lift coefficients
\begin{equation}
    C_D(t) = \frac{2}{U^2L}F_D(t),\quad C_L(t) = \frac{2}{U^2L}F_L(t)
\end{equation}
defined by the forces acting on the cylindrical obstacle
\begin{equation}
    \begin{aligned}
        F_D(t) & = \int_{\Gamma_\text{circle}} (-p(t)I+\nu\nabla\mathbf{v}(t))\cdot\mathbf{n}\cdot\mathbf{e}_1\;\mathrm{d}s\\
        F_L(t) & = \int_{\Gamma_\text{circle}} (-p(t)I+\nu\nabla\mathbf{v}(t))\cdot\mathbf{n}\cdot\mathbf{e}_2\;\mathrm{d}s
    \end{aligned}
\end{equation}
\subsubsection{Results}
We want to compare our results to the publicly available
results of the finite element software FEATFLOW~\citep{turek1996featflow}.
These are obtained with a Crank-Nicholson scheme with
time step size of \(k=1/1600\), which
corresponds to a Petrov-Galerkin discretization with \(cG(1)\) trial
and \(dG(0)\) test functions with \(M_{\texttt{FEAT}}=12800\)
temporal elements and DoFs. In space we compare to two
refinement levels, level 2 with \(2704\) DoFs and level 6 with
\(667,264\) DoFs.

In our own results we used \(M_\texttt{ideal}=256\) temporal elements
to see how a coarse solution with different temporal orders performs.
Specifically, we used \(r\in \{0,1,2\} \) with \( \{256,512,768\} \)
temporal DoFs. In space we used a fixed mesh with \(6100\) DoFs.

\begin{figure}[H]
    \begin{center}
        \includegraphics[page=1,width=0.8\textwidth]{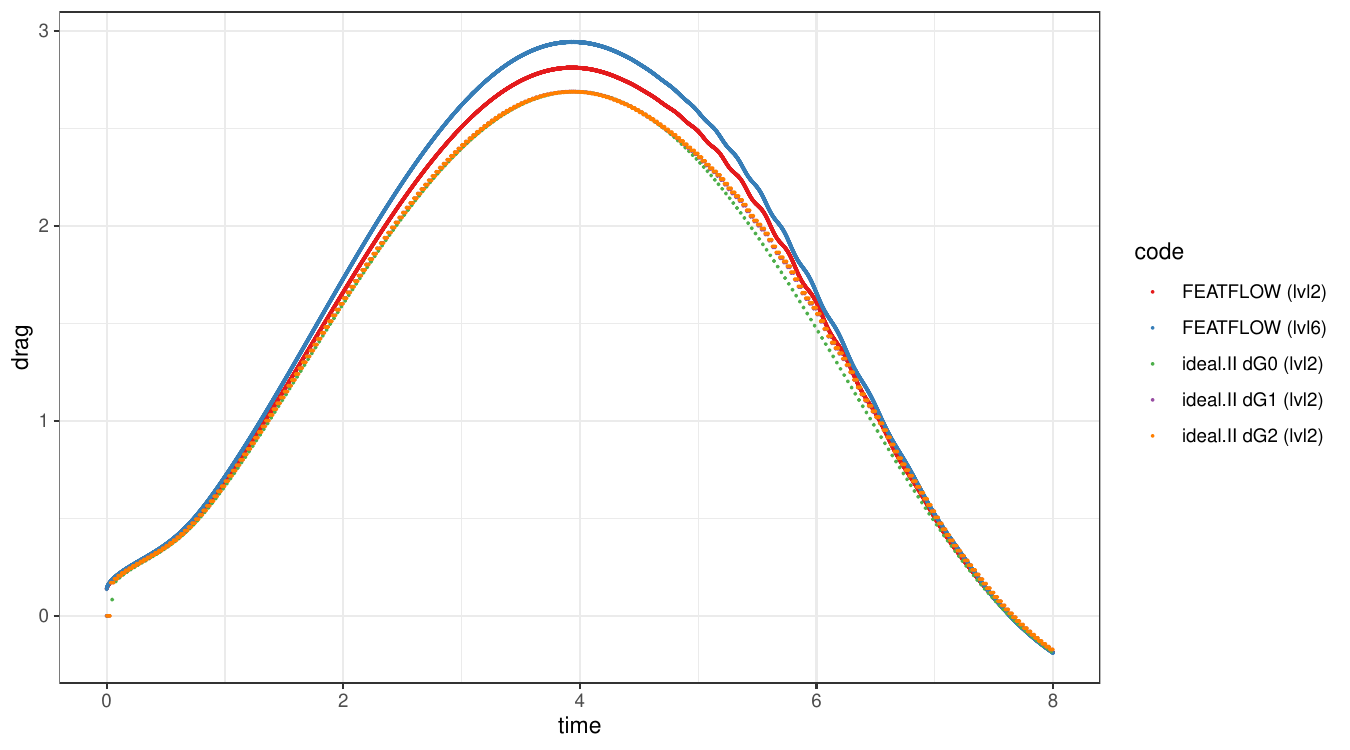}
    \end{center}
    \caption{Comparison of \(C_D(t)\) between \texttt{ideal.II} 
    with different \(dG(r)\) discretizations and \texttt{FEATFLOW} at different spatial
    refinement levels.}\label{fig:drag_results}
\end{figure}

\begin{figure}[H]
    \begin{center}
        \includegraphics[page=2,width=0.8\textwidth]{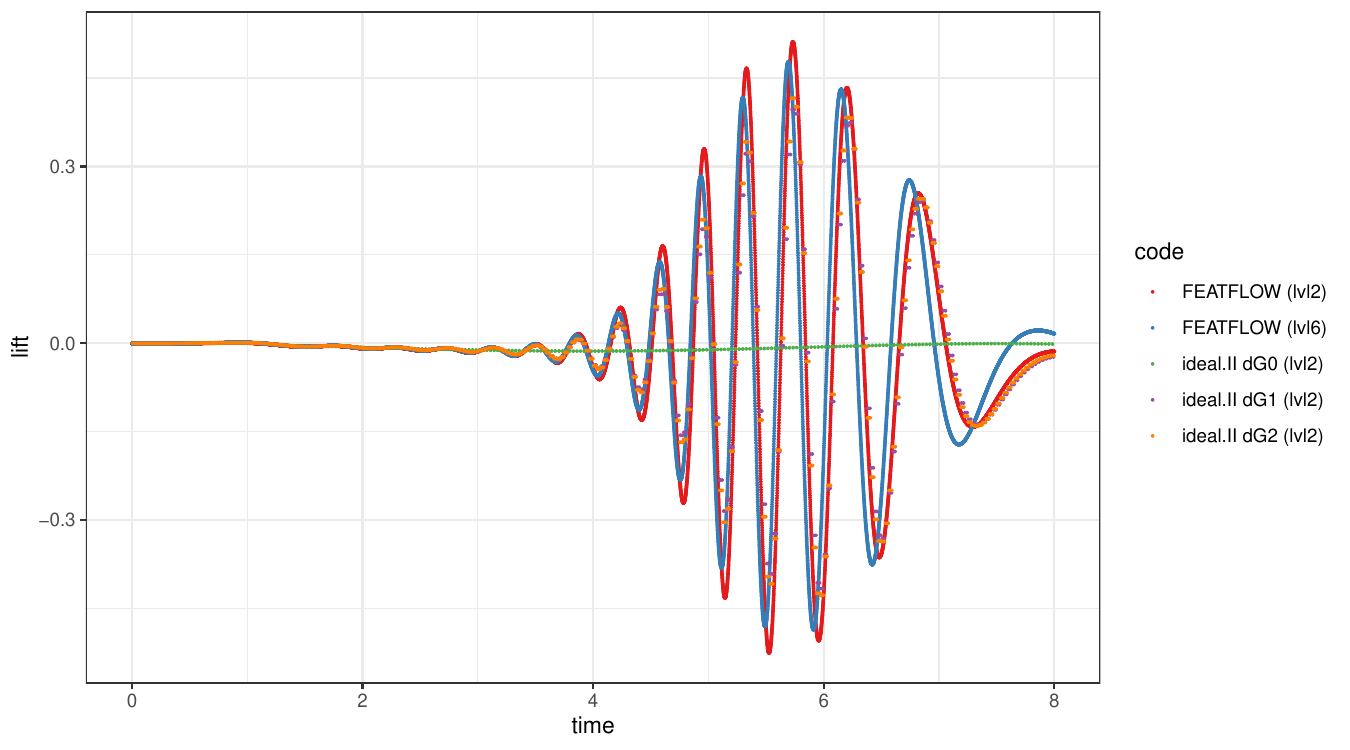}
    \end{center}
    \caption{Comparison of \(C_L(t)\) between \texttt{ideal.II} 
    with different \(dG(r)\) discretizations and \texttt{FEATFLOW} at different spatial
    refinement levels.}\label{fig:lift_results}
\end{figure}

Figure~\ref{fig:drag_results} shows that the drag coefficient
is mostly depending on the spatial discretization 
as the three \texttt{ideal.II} curves overlay each other 
until the \(dG(0)\) curve lowers faster after \(t\approx4.2\). 
In general, we also see that the slowing inflow 
in the second half introduces oscillations around the cylinder.
Only the \(dG(0)\) scheme, related to the backwards Euler scheme, 
produces a smooth curve due to numerical dampening.
For the lift coefficients in Figure~\ref{fig:lift_results} 
this dampening is even more pronounced.
For the other discretizations we see that the frequency of the 
oscillations depends on the spatial resolution.
For the amplitudes we see that, again, 
the peaks of the \texttt{ideal.II} results are below 
those of \texttt{FEATFLOW}.
Overall, we can conclude that the results are 
similar enough for \(r>0\) and that \texttt{ideal.II} 
is able to solve the benchmark,
even with the relatively low number of temporal DoFs.
\section{Outlook and further development}
In this section we will give an overview over the features we 
plan to implement into \texttt{ideal.II} in the future.

\subsection{Adaptive refinement and dynamic triangulation}
On the fixed triangulation, adaptive refinement is already possible
through adaptive refinement of the underlying meshed.
However, this interface is not convenient and 
local contributions on slabs need to be tallied manually
to obtain a single spatial set of error indicators 
and one indicator for each temporal element.
We want to provide an interface to make it easier to collect 
refinement information for fixed meshes.
Additionally, we plan to implement the dynamic triangulations 
as described in~\ref{subsec:disc}.

\subsection{PU DWR example step}
To show how to implement a problem with adaptive refinement, 
we plan to provide an example step that shows the application 
of the dual-weighted residual method with 
partition-of-unity localization (PU-DWR) 
for obtaining local error indicators.

This localization technique has originally been proposed 
in~\citep{RiWi15_dwr} for stationary problems 
and was extended to space-time settings in~\citep{ThiWi24}
and the PhD thesis of the first author~\citep{thiele2024}.

One additional component needed for this method are 
interpolation operators between different 
space-time finite element spaces. 
Since \texttt{deal.II} already offers such functions 
for the spatial and temporal components in the 
\texttt{FETools} namespace, these should be used 
internally.

\subsection{Further examples}
Currently, there are four examples:
\begin{enumerate}
    \item Introducing the basic concepts and solving the heat equation,
    \item Showing the solution of the coupled Stokes equations,
    \item Solving the heat equation in parallel via Trilinos and MPI,
    \item Solving a nonlinear problem, the Navier-Stokes equations.
\end{enumerate}
Apart from adaptivity, which will be shown in the aforementioned 
PU-DWR example step, there are other advanced topics and techniques
offered by \texttt{deal.II}, e.g.\ Multithreaded assembly using the
\texttt{WorkStream} class and automatic differentiation for nonlinear problems.
We could envision an example showing both for a nonlinear problem.

\subsection{Further linear algebra support}
In \texttt{deal.II} version 9.6 support for the newer Trilinos Tpetra
stack has been added. This stack not only supports MPI, 
but also multithreading and GPU calculations through Kokkos~\citep{Trott_Kokkos_3_Programming_2022}.

Additionally, \texttt{deal.II} offers an interface to the 
linear algebra and solvers provided by \texttt{PETSc}~\citep{petsc-efficient}.
However, to the best of the authors knowledge, 
the local indices of each MPI process have to be 
contiguous. This would require a reordering of the 
space-time indices such that all space-time indices belonging 
to the processor local portion of the spatial mesh 
fulfill this property.

\subsection{Automated testing and continuous integration}
Finally, we want to improve quality assurance of the library.

At the time of writing we have GitHub actions to check if the library can 
be build against the docker containers of a full *deal.II* installation 
in versions 9.3.0, 9.4.0, 9.4.2 and 9.5.0. 
Additionally, we check that all source code files are indented according to the 
clang-format specification checked into the repository.

We plan to add support for the GoogleTest framework to write unit 
and integration tests for the library. 

\begin{acks}
The author thanks Uwe K{\"o}cher (Hamburg),
Julian Roth (IfAM Hannover), Sebastian Kinnewig (IfAM), Thomas Wick (IfAM) 
and Christian Merdon (WIAS) for fruitful discussions
and Sebastian Kinnewig (IfAM) and Wolfgang Bangerth (Colorado State University) 
for contributions to the source code.
\end{acks}
\paragraph{Data availability}
The C++ library \texttt{ideal.II} described in this paper
is open source and can be found on 
\url{https://github.com/instatdealii/idealii}.
The presented simulation results have been
obtained with \texttt{ideal.II} version 1.0.0~\citep{idealii-1.0.0-zenodo}
built with 
\texttt{deal.II} version 9.5.2~cite{deal95}
and \texttt{Trilinos} version 15.1.0.
The results for the examples can be found on~\citep{idealii-step-3}
and~\citep{idealii-step-4} respectively.


\appendix
\bibliography{./lit}

\end{document}